\magnification=1200
%
%
\baselineskip=0.2in
\def\g{{\cal G}}

\def\Sp{{{\rm Sp}}}
\def\sp{{{\rm sp}}}
\def\so{{\rm so}}
\def\SO{{\rm SO}}
\def\tr{{\rm tr}}
\def\GL{{{\rm GL}}}
\def\I{{\cal I}}
\def\W{{\cal W}}
\def\IP{{{\cal IP}}}
\def\ltimes{\times'}
\def\H{{\cal H}}
\def\C{{\bf C}}
\def\R{{\bf R}}
\def\P{{\bf P}}
\def\Z{{\bf Z}}
\def\0{{\bf 0}}
\def\l{{\lambda}}
\def\a{\alpha}
\def\sym{{\rm sym}}
\def\Sym{{\rm Sym}}
\def\gl{{\rm gl}}
\def\lra{{\longrightarrow}}
\def\A{{\cal A}}
\def\Pf{{\rm Pf}}

\font\mysmall=cmr8 at 8pt

\centerline{{\bf ON THE DYNAMICS OF THE SELF GRAVITATING}}
\medskip
\centerline{{\bf ELLIPSOID IN $N$ DIMENSIONS AND}}
\medskip
\centerline{{\bf ITS DEFORMATION QUANTIZATION}}
\bigskip
\centerline{R. Fioresi\footnote*{Investigation supported by
the University of Bologna, funds for selected research topics.}}

\centerline{Dipartimento di Matematica, Universita' di Bologna}
\centerline{Piazza Porta San Donato 5, 40126 Bologna, Italy.}
\centerline{{\mysmall e-mail: fioresi@dm.unibo.it}}
\bigskip
{\bf 1. Introduction}
\medskip 
The self gravitating ellipsoid has been the subject of study by many
physicists and mathematicians. Newton first started the subject in the 
attempt to determine the eccentricity 
of the earth, which he modelled
using a rigidly rotating ellipsoid made of a fluid of uniform density and
subject only to its own gravity.
Later on McLaurin generalized and refined his calculation
of the eccentricity. It was the first time this number relative to
the shape of earth, was calculated with a reasonably accurate model.

\medskip

The dynamics of such an object was then studied by several people
among whom Dedekind, Riemann (often in the literature the self
gravitating ellipsoid is referred as the {\it Riemann
ellipsoid}) and more recently Chandrasekhar.
For an historical account of the development of this 
interesting subject see [Le], [Ch].

\medskip

Recently in [Ro] Rosensteel identified the phase space of the dynamical
system associated with the Riemann ellipsoid with the algebra
$\sym(3) \ltimes \gl(3)$, where $\sym(3)$ denotes the symmetric $3 \times 3$
matrices and $\gl(3)$ the $3 \times 3$ matrices.

\medskip

A generalization of some of his results was done by Carrero in [Ca].
Carrero introduced an arbitrary number of dimensions and 
studied in detail the coadjoint action of the group $G=\sym(n) \ltimes
\GL_+(n)$ on its Lie algebra $\sym(3) \ltimes \gl(3)$.
Global Darboux coordinates were explicitly calculated for any
coadjoint orbit of $G$. 

\medskip

Using such global
Darboux coordinates one can immediately write a deformation quantization
of the orbit using the Moyal Weyl type of deformation quantization. 
The existence of such a differential deformation, which is
unique up to gauge equivalence by Kontsevich's theorem [Kn],
does not however guarantee the existence of an algebraic deformation
quantization,
that is a deformation of the Poisson polynomial algebra of an orbit.
The explicit construction of such deformation will be the aim of this paper.

\medskip

This paper is organized as follows. 

\medskip

In \S 2 we study the coadjoint
action of $G=\sym(n) \ltimes \GL_+(n)$, reviewing some
of Carrero's results. We then show that there are no 
invariant polynomials with respect to this action. An
explicit description of the regular coadjoint orbits 
as algebraic varieties is given using semiinvariant polynomials.

\medskip

In \S 3 we prove the main result of this
paper, namely the existence of a deformation quantization of algebra of
polynomial functions of the regular coadjoint
orbits of $G$. The deformation is given explicitly
generalizing a construction
introduced in [FL1] in the case of a complex semisimple group.
The construction of the deformation is non trivial, since the method in [FL1]
depends in an essential way on the fact that the group is semisimple,
while our $G$ does not have this property.

\medskip

The existence of an algebraic deformation in this more
general setting 
suggests that a modification of the method in [FL1] could
possibly give quantization of more general Poisson algebraic variety.
We plan to explore this in a forthcoming paper.

\medskip

Acknoledgements. 
We want to thank Prof. Varadarajan for explaining the dynamics
of the Riemann ellipsoid.

\bigskip
{\bf 2. The coadjoint orbits of $G=\sym(n) \ltimes \GL_+(n)$}
\medskip
Let $G$ be the real Lie group $\sym(n) \ltimes \GL_+(n)$
with multiplication:
$$
(x,g)(y,h)=(x+g^{-1}y \check g,gh)
$$
where
$\sym(n)$ denotes the $n \times n$ real matrices, $\GL_+(n)$ the 
subgroup of $\GL(n)$ consisting of invertible
matrices with positive determinant and $\check g={g^t}^{-1}$. 

$G$ can be identified with a subgroup of $\Sp(n)$ in the following
way:
$$
\matrix{
G =\{(x,g)| x \in \sym(n), g \in \GL_+(n)\} \cong \cr\cr 
\{ \pmatrix{g & x{\check g} \cr 0 & {\check g}}
| x \in \sym(n), g \in \GL_+(n) \} \subset \Sp(n)
}
$$
It is immediate to check that with such identification we have that:
$$
\g =Lie(G) \cong \{ \pmatrix{a & b \cr 0 & -a^t} |
a \in \gl(n), b \in \sym(n)\} \subset \sp(n)
$$
\medskip

Let's consider the non degenerate form on $\sp(n)$: 
$$
<A,B>=\tr(AB).
$$
This form is still non degenerate on
$\g \times \g_- \subset \sp(n) \times \sp(n)$, where
$$
\matrix{
\g=\{\pmatrix{a & b \cr 0 & -a^t}| a \in \gl(n), b \in \sym(n) \} \qquad
\cr\cr
\g_-=\{\pmatrix{a & 0 \cr c & -a^t}| a \in \gl(n), c \in \sym(n) \} \qquad
}
$$
This allows us to identify ${\g}^* \cong \g_-$.

For brevity we will denote:
$$
\matrix{
\pmatrix{a & b \cr 0 & -a^t} \in {\g}^* \quad \hbox{with} \quad (b,a),
\qquad
\pmatrix{a & 0 \cr c & -a^t} \in {\g_-}^* \quad \hbox{with} \quad (c,a)
}
$$
and
$$
\pmatrix{g & x{\check g} \cr 0 & {\check g}} 
\in {G} \quad \hbox{with} \quad (x,g).
$$

\medskip
In the above notation we have that the adjoint and coadjoint actions
of $G$ on $\g$ and $\g^*$ respectively are given by: 
$$
\matrix{
Ad(x,g)(b,a)=({g}bg^{t}-\{gag^{-1}x+(gag^{-1}x)^t\}, gag^{-1})
\cr\cr
Ad^*(x,g)(c,a)=({\check g}cg^{-1},gag^{-1}+x{\check g}cg^{-1})
}
$$
\medskip
Define 
$$
\g^+=\{(c,a) \in \g^*| c \hbox{ positive definite}\} 
$$
We are now interested in the description of the coadjoint orbits
of $G$ of elements in $\g^+$. These are the orbits physically
interesting.
\medskip
Notice that 
$\g^+$ is an open set in $\g$ invariant under the coajoint action.
Let $O_{(c,a)}$ denote the coadjoint orbit of an element $(c,a)$.
Moreover one can immediately see that:
$$
O_{(c,a)}=O_{(I,d)}
$$
where $I$ is the identity matrix and $d \in \gl(n)$.
\medskip
{\bf Lemma (2.1)}. 
{\it
Let $(c,a) \in \g^+$ and let $O_{(c,a)}$ be the coadjoint 
orbit of $(c,a)$.

1)  If $n=2k+1$, there exists a unique element $(I,H)$ such
that $O_{(c,a)}=O_{(I,H)}$ with
$$
H= \pmatrix{ 0   & -\l_1 & \dots & & & &\cr
           -\l_1 &   0   &       & & & &\cr
                 & &               & & & &\cr   
                 & & \vdots   & & & &\cr
                 & & & & & &\cr
                 & & &   & 0 & \l_k & \cr  
                   & &   & & -\l_k & 0 & \cr 
                 & & & & & & 0  }   
$$
and $\l_1 \geq \l_2 \geq \dots \geq \l_k$.

2) If $n=2k$ there exists a unique element $(I,H)$ such
that $O_{(c,a)}=O_{(I,H)}$ with
$$
H= \pmatrix{ 0   & -\l_1 & \dots & & & &\cr
           -\l_1 &   0   &       & & & &\cr
                 & &               & & & &\cr
                 & & \vdots   & & & &\cr
                 & & & & & &\cr
                 & & &   & 0 & \l_k & \cr
                   & &   & & -\l_k & 0 & }
$$
and $\l_1 \geq \l_2 \geq \dots \geq l_k$.

3) If $\l_1 > \l_2 > \dots > \l_k$, 
$$
O_{(c,a)}=G/\SO(2) \times \dots \times \SO(2)
$$
hence dim $O_{(c,a)}=n^2-k$.
}

{\bf Proof.} 
See [Ca].
\medskip

Let $\H$ be the Cartan subalgebra of $\so(n)$ defined in Lemma (2.1).
Let $\I(\H)$ be the algebra of invariant polynomials under the action
of the Weyl group $\W$ of $\so(n)$. We know that this algebra is the
same as the algebra of polynomials on $\so(n)$ invariant under
the adjoint action. Since every $G$ orbit in $\g^+$ meets $\H$
in a $\W$ orbit we have that an invariant function on $\g^+$ is determined
uniquely by its restriction to $\H$. Let $\I(\g^+)$ be the algebra
of invariant functions on $\g^+$.
We have the following result.

\medskip
{\bf Theorem (2.2)}. 
{\it 
1) If $n=2k$ then $\I(\H)$ is generated by:
$$
\matrix{
\tilde f_i=\sum_{j=1}^k \l_j^{2i} \quad i=1 \dots k-1, \quad
\tilde f_k=\l_1 \dots \l_k
}
$$
and $\I(\g^+)$ is generated by:
$$
\matrix{
f_i={1 \over 2^{2i}}\tr(cac^{-1}-a^t)^{2i}, \quad
i=1 \dots k-1, \quad \Pf=\pi(ac^{-1})det(c)^{1/2} 
\cr\cr
}
$$
with $\pi(A)$ denoting the Pfaffian of the matrix $1/2(A-A^t)$.

2) If $n=2k+1$, then 
$\I(\H)$ is generated by:
$$
\matrix{
\tilde f_i=\sum_{j=1}^k \l_j^{2i} \quad i=1 \dots k
}
$$
and $\I(\g^+)$ is generated by:
$$
\matrix{
f_i={1 \over 2^{2i}}\tr(cac^{-1}-a^t)^{2i}, \quad 
i=1 \dots k. \cr\cr
\cr\cr
}
$$
}

{\bf Proof}. See [Ca]. 

\medskip

Notice that there are rational and irrational functions in 
$\I(\g^+)$ that are not polynomial in $(c,a)$.

\medskip
We now want to determine the subring of invariant polynomials
$\IP(\g^+) \subset \I(\g^+)$.

\medskip
{\bf Proposition (2.3)}. 
{\it If $n=2k$, the ring $\IP(\g^+)$ of invariant polynomials 
on ${\g^+}^*$ consists only of constants.}

{\bf Proof}. 
We first observe that 
$$
\I(\g^+)=span_{l \in \Z} \tr(cac^{-1}-a^t)^{2l}
$$
Assume that the polynomial $f(a,c) \in  \I(\g^+)$. This means that
$$
f(c,a)=\sum_{\l_l \in \C} \l_l \tr(cac^{-1}-a^t)^{2l}
$$
Let $-2M$ be the highest negative degree of $det(c)$. Then we
have the equation between polynomials:
$$
det(c)^{2M}f(c,a)=\sum_{\l_l \in \C} \l_l det(c)^{2M-2l}
\tr(caC^t-a^t)^{2l}
$$
where $C$ is the matrix of the algebraic complements.

Observe that an invariant polynomial must depend on both
$a$ and $c$. We will prove that it depends only on $a$ reaching
a contradiction. Let $\deg_{ij}$ denote the degree in $c_{ij}$ of
a generic polynomial, where $c_{ij}$ is the $(i,j)$ entry of
the matrix $c$.
\medskip
{\bf Claim}. $\deg_{ij}(\tr(caC^t-a^t)^{2l} \leq 2l$.

Given a matrix $x$ let's associate to it another matrix $(x)_{\deg_{ij}}$
whose $(k,l)$ entry is $\deg_{ij}(x_{kl})$. 

One can easily see the following:
$$
\matrix{
(c)_{\deg_{ij}}=E_{ij} \cr\cr
(a)_{\deg_{ij}}=\0 \cr\cr
(C^t)_{\deg_{ij}}=\sum_{1 \leq r,s \leq k, r\neq i,s \neq j} E_{rs}
}
$$
where $E_{ij}$ denotes the elementary $k \times k$ matrix having
1 in the $(i,j)$ position and 0 everywhere else and $\0$ denotes the
$k \times k$ zero matrix.

By induction one gets:
$$
((caC^t)^m)_{\deg_{ij}}=
\sum_{1 \leq r,s \leq k, r\neq i,s \neq j} mE_{rs}  
+ \sum_{1 \leq r \leq k, r \neq i} (m-1)E_{rj}+
\sum_{1 \leq s \leq k,s \neq j} (m+1)E_{rs} 
+ mE_{ij}
$$
From which we have $\deg_{ij}(caC^t)^{2l} \leq 2l+1$, $i \neq j$, 
$\deg_{ii}(caC^t)^{2l} \leq 2l$.

Now we compute $\deg_{ij}(f(c,a))$. By the claim we have
$$
\deg_{ij}(\sum_{\l_l \in \C} \l_l det(c)^{2M-2l}
\tr(caC^t-a^t)^{2l}) \leq \cases{2M+1 & if $i \neq j$ \cr
				 2M & if $i=j$}.
$$
But observe that
$$
\deg_{ij}(det(c)^{2M})=2M.
$$
Hence $\deg_{ij}(f(c,a)) \leq 1$ if $i \neq j$,
$\deg_{ii}(f(c,a))=0$. This implies that $f(c,a)$ does not
depend on $c_{ii}$. Now assume
$$
f(c,a)=\sum_{b_{ij} \in \C[a], i \neq j} b_{ij}c_{ij}.
$$
where $\C[a]$ denotes the ring of polynomials in the indeterminates $a_{ij}$'s.
Now choose $(x,g) \in G$ such that $\check g_{j_0i_0}\neq 0$,
$g_{j_0j_0}^{-1} \neq 0$.
It is a simple computation to check that
$$
(\check g c \check g^{-1})_{j_0j_0} \neq 0.
$$
This implies that $f(c,a)$ must also depend on $c_{j_0j_0}$ which is
a contradiction. So $f(c,a)$ must depend only on $a$, but this is
not possible since it is invariant, unless it is a constant.
QED.

\medskip
We now want to ask whether there are semiinvariants for the coadjoint
action.

Define the polynomials:

Odd case:
$$
h_i=det(c)^{2i}\tr(cac^{-1}-a^t)^{2i} \quad i=1 \dots k
$$
Even case:
$$
\matrix{
h_i=det(c)^{2i}\tr(cac^{-1}-a^t)^{2i} \quad i=1 \dots k-1, \cr\cr
h_k=det(c)^2 \pi(ac^{-1})^2 det(c)
}
$$
where $\pi(ac^{-1})=\Pf(1/2(ac^{-1}-(ac^{-1})^t)$ and $\Pf$
denotes the Pfaffian.

One can easily check that these algebraic functions on $\g^*$:
are semiinvariant for the coadjoint action. In fact:
$$
Ad^*_{(x,g)}h_m=det(g)^{-4m} h_m
$$
It would be interesting to determine all semiinvariants.
\medskip
We now would like to describe the coadjoint orbits as algebraic varieties.
For this reason is now more convenient to look at their complexification.
Let $O_{(c,a)}^\C$ denote the complexification of the orbit $O_{(c,a)}$.

It is clear that we cannot describe the ideal of the orbit
 in the same way as in the semisimple
case. In fact in that case we have that the ideal of a given regular
orbit is simply given by the polynomials that Chevalley generators of the
ring of invariant polynomials equal to constant ([Ko]). We will need to
use the semiinvariant polynomials.
\medskip

{\bf Theorem (2.4).}
{\it Given a regular orbit $O_{(c,a)}$, its ideal is given by:

1) If $n=2k+1$, $(h_1-\a_1det(c)^{2}, \dots h_k-\a_kdet(c)^{2k})$

2) If $n=2k$, $(h_1-\a_1det(c)^{2}, \dots h_k-\a_kdet(c)^{2})$

where $\a_i$ for $i=1 \dots k-1$
is the constant value of the rational function
$\tr(cac^{-1}-a^t)^{2i}$ on the orbit.
$\a_k$ is the the constant value of the rational function
$\tr(cac^{-1}-a^t)^{2k}$ 
on the orbit if $n=2k+1$, while if $n=2k$ it is the
constant value of $\pi(ac^{-1})^2 det(c)$ on the orbit.}

{\bf Proof.} Let $r_1 \dots r_k$ be the generators of the
ideal.
Since the orbit is a non singular algebraic variety
of dimension $n^2-k$
it is enough to prove that the differentials 
$dr_1 \dots dr_k$
are linearly independent over every point of $O_{(c,a)}$.
Since there is a diffeomorphism that brings any point of 
$O_{(c,a)}$ into any other point, it is enough to prove the
differentials are linearly independent over points in $\H$.

Observe that $d(r_1|_\H), \dots d(r_k|_\H)$ are linearly
independent over every regular point of $\H$ ([Va2]).

Hence it is simple to see that also $(dr_1)|_\H \dots (dr_k)|_\H$
are linearly independent over every regular point of
$\H$. Q.E.D.
\bigskip
{\bf 3. Deformation quantization of regular coadjoint orbits} 

\medskip
We would like to construct a deformation quantization of the
algebra of regular functions on a regular coadjoint orbit of a
Lie group under certain hypothesis listed below, which are satisfied by 
$G=\sym(n) \ltimes \GL_+(n)$. Our construction is a generalization of 
the one described in [FL1] where $G$ was assumed to
be complex semisimple.

Let's recall the basic definitions.
\medskip
{\bf Definition (3.1)}. Given a real (or complex) Poisson algebra 
$\P$, a {\it formal deformation} or a {\it deformation quantization} of $\P$
is an associative algebra $\P_h$ over $\R[[h]]$ (or over \C[[h]]), 
where $h$ is a formal
parameter, with the following properties:

\noindent a. $\P_h$ is isomorphic to $\P[[h]]$ as a  $\R[[h]]$-module
(or as a $\C[[h]]$-module).

\noindent b. The multiplication $*_h$ in $\P_h$ reduces  mod($h$) to 
the one in $\P$.

\noindent c. $\tilde F *_h \tilde G -\tilde  G *_h \tilde F = h\{F,G\}$
mod $(h^2)$, where $\tilde F, \tilde G \in \P_h$ reduce to $F,G \in \P$ 
mod($h$) and $\{\, ,\,\}$ is the Poisson bracket in $\P$.

\medskip
This definition makes sense also if
we substitute $\C[[h]]$ by $\C[h]$. In this case we 
will say that $P_h$ is a {\it $\C[h]$-deformation}.

Notice that a $\C[h]$-deformation extends immediately to a formal one
by tensoring by $\C[[h]]$, 
but the converse is not always true. 
Moreover $\C[h]$-deformation can be specialized to any value of the
parameter $h$.

\medskip
Let's now make the following assumptions.

1. $\g$ is a finite dimensional complex Lie algebra of a complex Lie group $G$.

2. $p_1 \dots p_m \in \C[\g^*]$ are semiinvariant polynomials with
respect to the coadjoint action.

3. $dp_1 \dots dp_m$ are linearly independent over points where
$p_1=\dots=p_m=0$.

4. The set of zeros of $p_1 \dots p_m$ is an algebraic Poisson variety
with bracket induced by the one in $\g^*$.
\medskip
In these hypothesis we will construct a deformation quantization of the
algebraic Poisson variety described by the ideal $(p_1 \dots p_m)\subset 
\C[\g^*]$, i.e. a formal deformation of the Poisson algebra
$\C[\g^*]/(p_1 \dots p_m)$.

\medskip
{\bf Observation (3.2)}.
Notice that if we take
$\g=Lie(\sym(n) \ltimes \GL_+(n))$, $m=k$ and
$p_i=h_i-\a_idet(c)^{2i}$,
$i=1 \dots k-1$, and if $n=2k+1$
$p_k=h_k-\a_idet(c)^{2k}$, if $n=2k$ $p_k=h_k-\a_idet(c)^{2}$
(for the notation see \S 2) the hypothesis
(1), (2), (3), (4) listed above are satisfied.
Hence the procedure described below will give us a deformation
quantization of the algebra of polynomial function on a
regular coadjoint orbit of $\sym(n) \ltimes \GL_+(n)$,
$O^\C_{(c,a)}$, set of zeros of such $p_1 \dots p_k$.
\medskip
Let's denote by $T_A(V)$ the full tensor algebra of
a complex vector space $V$ over a $\C$-algebra $A$. Let $\g=Lie(G)$. 
Consider the  proper two sided ideal in $T_{\C[h]}(\g)$
$$
{\cal L}_h=\sum_{X,Y \in \g}
T_{\C[h]}(\g) \otimes(X \otimes Y - Y \otimes X - h[X,Y])
\otimes T_{\C[h]}(\g)
$$
\medskip
$U_h$ is a free $\C[[h]]$-module, in particular it is torsion free
([FL1] Proposition (3.2)).
\medskip
Define $U_h=_{def} T_{\C[h]}(\g)/\l_h$. $U_h$ is the 
universal enveloping algebra of the Lie algebra $\g_h=\C[h] \otimes_{\C} \g$
with Lie bracket 
$$
[p(h) X,q(h) Y]_{h}=p(h)q(h)[X,Y]
$$

Let $\C[X]$ denote the regular (polynomial) functions on a variety $X$.
Define $\C_h[\g^*]=\C[h] \otimes \C[\g^*]$. Observe that
$\C_h[\g^*] \cong \C[\g_h^*]$.
Let $\Sym$ denote the symmetrizer map, 
$\Sym:\C_h[\g^*] \lra U_h$ ([Va1] pg. 180).

Denote also with $I_h$ the two-sided ideal in $U_h$ generated by
$P_1=\Sym(p_1), \dots P_m=\Sym(p_m)$.

\medskip
{\bf Proposition (3.2)}. {\it If $p \in \C[\g^*]$ is semiinvariant
with respect to the coadjoint action i.e.:
$Ad^*(g)p=f(g)p$ for all $g \in G$,
where $f$ is a function depending only on $g$,
then 
$$
[X,\Sym(p)]=F(X) \Sym(p) \qquad 
\hbox{for all } X \in U_h,
$$  
with $F$ scalar function depending only on $X$.}

{\bf Proof}. Direct calculation.
\medskip
{\bf Lemma (3.3)}. 
{\it
Let $r$  be a fixed positive integer and let all the notation be as above.
Let
$$
\sum_{1 \leq i_1\leq  \cdots \leq i_r \leq m} 
a_{i_1\dots i_r}p_{i_1}
\dots p_{i_r}=0
$$
with $ a_{i_1\dots i_r}\in \C[\g^*]$.
Then $a_{i_1\dots i_r} \in (p_1,\dots,p_m) \subset \C[\g^*]$.}

{\bf Proof}. The proof is the same as in Proposition (3.8), [FL1].

\medskip
{\bf Lemma (3.4)}.
{\it Let $k$ be a fixed integer and let
$$
\sum_{i_1\leq\cdots i_k\leq m} A_{i_1\dots
i_k}P_{i_1}\cdots P_{i_k}
\equiv 0 \qquad \hbox{mod} h
$$
where $A_{i_1\dots i_k} \in U_h$ and $P_i=\Sym(p_i)$.  

Then
$$
\matrix{
\sum_{i_1\leq\cdots i_k\leq m} A_{i_1\dots
i_k}P_{i_1}\cdots P_{i_k}
=h\sum_{i_1\leq\cdots i_k\leq m} B_{j_1 \dots j_l,i_1\dots
i_k} 
\cr \cr
P_{j_1}\cdots P_{j_l}
P_{i_1}\cdots P_{i_k}
}
$$
}

{\bf Proof}. This is the same as Lemma (3.9), [FL1].

\medskip
{\bf Lemma (3.5)}.
{\it If $hF \in I_h$ then $F \in I_h$. In other words, 
$U_h/I_h$ is torsion free.}

{\bf Proof}.
Since $hF \in I_h$ and since Proposition (3.2) we
can write:
$$
hF = \sum A_i P_i
$$
We have $\sum A_i P_i \equiv 0$ mod$h$. Hence, by Lemma (3.4)
and also by the fact that $U_h$ is torsion free
we have our result.

\medskip

We now want to construct a basis for the torsion free $\C[h]$-module $U_h/I_h$.

\medskip
Let's fix a basis $\{X_1,\dots, X_n\}$ of $\g$ and let $x_1,\dots , x_n$ 
be the corresponding elements in $\C[\g^*]$. With this choice $\C[\g^*]
\cong \C[x_1,\dots , x_n]$.
Let $\{x_{i_1},
\dots , x_{i_k}\}_{(i_1,\dots , i_k) \in \A}$ be  a basis in
of $\C[\g^*]/(p_1 \dots p_m)$ as $\C$-module, where $\A$
is a set of multiindices
appropriate to describe the basis. In particular, we can take them such 
that $i_1 \leq \cdots\leq i_k$.

\medskip

{\bf Proposition (3.6)}.
{\it The monomials $\{X_{i_1}\cdots X_{i_k}\}_{(i_1, \dots , i_k) \in \A}$
are a basis for  $U_h/I_h$.}

{\bf Proof}. The proof is exactly the same as the one in
Proposition (3.11) and (3.13) in [FL1].
\medskip

{\bf Theorem (3.7)}. {\it Let the notation be as above. 

\noindent
1. $U_h/I_h$ is a $\C[h]$-deformation of
$\C[\g^*]/(p_1 \dots p_m)$.

\noindent 
2. $U_h/I_h \otimes \C[[h]]$ is a deformation quantization
of $\C[\g^*]/(p_1 \dots p_m)$. 
}

{\bf Proof.} Immediate from previous lemmas.

\medskip 
\bigskip
{\bf References}
\medskip
\item{[Ca]} 
J. Carrero  {\it Lie aspects of the dynamics of the self-gravitating
ellipsoid in $n$-dimensions}, JMP, {\bf 42}, no. 4, 1761-1778, 2001.
\vskip1ex
\item{[Ch]} 
S. Chandrasekhar  {\it Ellipsoidal figures of equilibrium}, New York
Dover publications, 1987.
\vskip1ex
\item{[FL1]} R. Fioresi  and M. A. Lled\'o  {\it On the deformation
Quantization of Coadjoint Orbits of Semisimple Lie Groups},
Pacific J. of Math.,{\bf 198}, No. 2, 411-436, 2001.
\vskip1ex
\item{[FL2]} 
R. Fioresi  and M. A. Lled\'o  {\it A comparison between star products
on regular orbits of compact Lie groups}, J. Phys. A: Math. Gen., {\bf 35},
1-13, 2002. 
 \vskip1ex
\item{[FLL]}  
R. Fioresi, A. Levrero, M.A. Lledo {\it Algebraic and differential
       star products on regular orbits of compact Lie groups}, 
       preprint, 2001.
\vskip1ex
\item{[Kn]} M. Kontsevich {\it Deformation Quantization of Poisson Manifolds}.
Preprint, \break q-alg/9709040, (1997).

\item{[Ko]} B. Kostant  {\it Lie group representations on polynomial
rings}. Am. J. Math. {\bf 85} 327, 1978.
\vskip1ex
\item{[Le]} 
N. Lebovitz  {\it The mathematical development of the classical ellipsoid},
Int. J. Eng. Sci., {\bf 36}, 1407-1420, 1998.
\vskip1ex
\item{[Ro]} 
G. Rosensteel  {\it Rapidly rotating nuclei as Riemann ellipsoids}, 
Ann. Phys. (N.Y.), {\bf 86}, 230-291, 1988.
\vskip1ex
\item{[Va1]} V. S. Varadarajan  {\it Lie groups, Lie algebras and their
representations}, GTM, Springer Verlag, 1974.
\vskip1ex
\item{[Va2]} V. S. Varadarajan  {\it On the ring of invariant polynomials 
on a semisimple Lie algebra}. Amer. J. Math, {\bf 90}, 1968.

\end